\documentclass[usletter]{article}

\usepackage{amsmath, amssymb}

\newtheorem{theorem}{Theorem}
\newtheorem{corollary}[theorem]{Corollary}

\newtheorem{proposition}[theorem]{Proposition}

\newcommand{\sector}[1]{
  \refstepcounter{section}
  \setcounter{equation}{0}
  \setcounter{theorem}{0}
  \mbox{\ }\smallskip
  {\begin{center}
    \bf \thesection.\ {#1}
  \end{center}}}

\newcommand{\eofproof}{{ \hspace{0.2cm}\\
        \mbox{\,}\hfill\framebox[0.3cm]{\rule{0cm}{.1cm} } }}
\newcommand{\Pf}[1]{{\noindent \bf Proof{#1}\ }}

\newcommand{\eqdef}{{\ \stackrel{\mathrm{def}}{=}\ }}

\begin{document}

\title{Hierarchical size-structured populations: The linearized semigroup approach}
\author{J\'{o}zsef Z. Farkas\thanks{E-mail: jzf@maths.stir.ac.uk}\ {\normalsize $^1$}\mbox{\ }
and Thomas Hagen\thanks{E-mail: thagen@memphis.edu}\ {\normalsize $^2$} \mbox{\
}
    \mbox{\ }\\
    \mbox{\ }\\
    {\normalsize $^1$\,Department of Computing Science and Mathematics,
     University of Stirling,}\\
    {\normalsize  Stirling, FK9 4LA, UK}  \\
{\normalsize $^2$\,Department of Mathematical Sciences, The University of Memphis,}\\
    {\normalsize  Memphis, TN 38152, USA} }
\date{September 2008}

\maketitle

\begin{abstract}
In the present paper we analyze the linear stability of a
hierarchical size-structured population model where the vital
rates (mortality, fertility and growth rate) depend both on size
and a general functional of the population density (``environment"). We derive
regularity properties of the governing linear semigroup, implying
that linear stability is governed by a dominant real eigenvalue of
the semigroup generator, which arises as a zero of an associated
characteristic function. In the special  case where neither the
growth rate nor the mortality depend on the environment, we
explicitly calculate the characteristic function and use it to
formulate simple conditions for the linear stability of population
equilibria. In the general case we derive a dissipativity
condition for the linear semigroup, thereby characterizing
exponential stability of the steady state.

\end{abstract}

\sector{Introduction}

In the last three decades nonlinear age- and size-structured
population models have attracted a lot of interest both among
theoretical biologists and applied mathematicians. Traditionally,
structured population models have been formulated as partial
differential equations for population densities. Starting with the
seminal work \cite{GM}, researchers have been developing
and analyzing various physiologically structured population
models. We refer here to the well-known monographs
\cite{CUS,I,MD,WEB}.

Diekmann et al. have been developing a general mathematical
framework for modeling structured populations, see for example
\cite{DGH,DGM}. One of their most important
recent results is that the qualitative behavior of nonlinear
physiologically structured population models can be studied by
means of linearization \cite{D1,D2}. In other words, they have
proven for a very general class of physiologically structured
population models that the nonlinear stability/instability of a
population equilibrium is completely determined by its linear
stability/instability. Such a fundamental result, often referred
to as ``the Principle of Linearized Stability", has been shown
previously for several concrete age- and size-structured models
\cite{GM,K,PR1,TZ,WEB}.

Following the lead of \cite{PR1} and \cite{WEB}, we
successfully applied linear semigroup methods to formulate
biologically interpretable conditions for the linear
stability/in\-sta\-bi\-li\-ty of equilibria of several structured
population models \cite{FH1,FH3,FH2}. In these problems the vital
rates depend on size or age and on the total population size, in
general. Hence it is assumed that any effect of intraspecific
competition on individual behavior is primarily due to a change in
population size and every individual in the population can
influence the vital rates of other individuals, a scenario
commonly referred to as ``scramble competition".

In other scenarios competition among individuals is based upon
some hierarchy in the population which is often related to the
size of individuals. In this case the nonlinearity (environmental
feedback) in the model is incorporated through infinite
dimensional interaction variables. A simple example for this
situation is given by a forest consisting of tree individuals in
which the height of a tree determines its rank in the population
\cite{KR}. Taller individuals have higher efficiency when
competing for resources such as light, while individuals of lower
rank cannot affect the vital rates of individuals of higher rank.
This scenario is the so called ``contest competition". Both
discrete time  and continuous hierarchical structured models have
been developed, see \cite{CJ} and the references therein.

Of interest in this work is the stability analysis of population
equilibria  by means of linearization of a continuous quasilinear
size-structured model, recently discussed in \cite{ADH}. In
this model the density evolution of individuals of size $s$ is
assumed to be governed by the following quasilinear partial
differential equation
\begin{equation}
u_t(s,t)+(\gamma(s,Q(s,t))\,u(s,t))_s+\mu(s,Q(s,t))\,u(s,t)=0,\label{pde}
\end{equation}
defined for $0\le s \leq m<\infty$ and $t>0$. The density of zero
(or minimal) size individuals is given by the nonlocal boundary
condition
\begin{equation}
u(0,t)=\int_0^m \beta(s,Q(s,t))\,u(s,t)\,ds,\quad t>0.
\label{boundary}
\end{equation}
The quantity $m$ denotes the maximum size of individuals.
The initial condition takes the form
\begin{equation}
u(s,0)=u_0(s),\quad   s\in [0,m]. \label{initial}
\end{equation}
Here $\beta,$ $\mu$ and $\gamma$ denote the fertility, mortality
and growth rate of individuals, respectively. We assume that these
vital rates depend on the individual size $s$ and on the
environment
\begin{equation}
Q(s,t)=\alpha\,\int_0^s w(\eta)\,u(\eta,t)\,d\eta+\int_s^m
w(\eta)\,u(\eta,t)\,d\eta,\quad 0\leq s\leq m,\quad t\geq
0.\label{density}
\end{equation}
The constant $\alpha$ is a parameter in $[0,1]$ measuring the
degree of hierarchy in the population, while the function $w$
represents a positive weight. For example in case of a tree
population where taller individuals overshadow smaller individuals
\cite{KR} the vital rates of an individual of size $s$ are
reasonably assumed to depend on the cumulative leaf area of
individuals of size $s$ or larger, modeled by the function
\begin{equation}
Q(s,t)=\int_s^m w(\eta)\,u(\eta,t)d\,\eta.
\end{equation}
Here $w$ is an appropriately chosen weight function. Hence in this
situation the parameter $\alpha$ would be $0$. The case $\alpha=1$
(which represents scramble competition) has been treated in
 detail in  \cite{FH1}.

We impose the following regularity conditions on the model ingredients:
\begin{itemize}
\item $\mu=\mu(s,Q)\in C([0,m];C^1[0,\infty))$, $\mu\geq 0$
\item $\gamma=\gamma(s,Q)\in C^1([0,m]; C^1[0,\infty))\cap
C([0,m]; C^2[0,\infty))$, $\gamma>0$
\item $\beta=\beta(s,Q)\in C([0,m];C^1[0,\infty))$, $\beta\geq 0$
\item $w=w(s)\in C^1([0,m])$, $w> 0$.
\end{itemize}
These assumptions are tailored toward the linear analysis of this
work. They might, however, not suffice to guarantee the existence
and uniqueness of solutions of Eqs.~\eqref{pde}--\eqref{density}.
Well-posedness of structured partial differential equation models
with infinite dimensional environmental feedback variables is in
general an open question, although conditions for the global
existence of weak solutions in the case discussed here are given
in \cite{ADH}. It has recently been shown \cite{AI,KR} that
the population model \eqref{pde}--\eqref{density} may exhibit a
more complicated dynamical behavior than the simple
size-structured model of scramble competition. In particular, in
\cite{AI} it was demonstrated both analytically and
numerically that a singular solution of
\eqref{pde}--\eqref{density} containing a Dirac delta mass
component can emerge if the growth rate $\gamma$ is not a
decreasing function of the environment $Q$.

For a more realistic description of real populations in a specific
setting, one will have to modify the assumptions on the vital
rates above. For example, one would possibly demand that
$\displaystyle{\lim_{s\rightarrow m}\mu(s,.)=\infty}$, thus
modeling a gradual rather than instantaneous reduction in the
numbers of individuals reaching maximum size $m$.

The size-structured model \eqref{pde}-\eqref{density}
is often considered (see \cite{AI,ADH}) with a boundary condition of the
form
\begin{align}
\gamma(0,Q(0,t))\,u(0,t)=C(t)+\int_0^m
\beta(s,Q(s,t))\,u(s,t)\,ds,\quad t>0.\label{oldboundary}
\end{align}
In \eqref{boundary} we have taken $C\equiv 0$ and incorporated the
growth rate $\gamma(0,Q(0,t))$ on the left of \eqref{oldboundary}
in the birth rate $\beta(s,Q(s,t))$ on the right of
\eqref{boundary}, assuming that zero size individuals grow
instantaneously. This assumption seems reasonable, for example in
case of a forest population. It is then clear that the two
boundary conditions are equivalent (in the case $C\equiv 0$
treated here).  We have observed, however, that \eqref{boundary}
is better suited for analytical work \cite{FH1,FH3,FH2}. As recent
results indicate \cite{FJ4,FH2}, the introduction of a positive
inflow $C$ may have a significant influence on the linearized
dynamical behavior of \eqref{pde}--\eqref{density}. A
comprehensive study of the effects of a positive inflow in
hierarchical populations is left for future work.

The study of hierarchical models in the literature  is largely
based on a decoupling of the total population quantity from the
governing equations and a transformation  of the nonlocal partial
differential equation \eqref{pde} into a local one
\cite{CS2,C1,KR}. This technique allows to prove well-posedness
and to study the asymptotic behavior of solutions by means of ODE
methods. For Eqs.~\eqref{pde}--\eqref{density} this transformation
fails since the vital rates depend on both size $s$ and on the
environment $Q(s,t)$. Therefore it seems unavoidable to study the
original partial differential equation \eqref{pde} with the
nonlocal integral boundary condition \eqref{boundary} directly.
This approach is based on a linearization of the governing
equations about steady state  \cite{FJ4,FH1,FH3,FH2,PR1,WEB}.
While Sections~3 through 5 exploit spectral theoretic and
structural properties of the governing linear semigroup extending
related results in \cite{FH1,FH3,FH2}, Section 6 gives a new characterization of
asymptotic stability of the semigroup in terms of a dissipativity
criterion. This idea was previously introduced and employed in
\cite{Ha2,Ha4} for elongational flow
problems.

\sector{The linearized system}

Eqs.~\eqref{pde}--\eqref{density} have obviously the trivial
solution $u_*\equiv 0$. Realistically we also expect additional
positive (continuously differentiable) solutions $u_*>0$. In the following we formulate a necessary
condition  for the existence of a positive equilibrium solution of
problem \eqref{pde}-\eqref{density}.

\begin{proposition}\label{existstatsol}
If $u_*$ is a positive stationary solution of problem
\eqref{pde}--\eqref{density}, then the function
 $Q_*$, defined by
 \begin{equation}\label{Qdef}
Q_*(s)   =\alpha\,\int_0^s w(\eta)\,u_*(\eta)\,d\eta+\int_s^m
w(\eta)\,u_*(\eta)\,d\eta,
\end{equation}
 satisfies  the equation
\begin{equation}\label{netrep}
R(Q_*)=1,
\end{equation}
where $R:C([0,m])\to\mathbb{R}$ is the inherent net reproduction
rate
\begin{equation}
R(Q) \eqdef\int_0^m\beta(s,Q(s))\,\pi(s,Q)\,ds\label{netrep2}
\end{equation}
and the operator $\pi$ is given for $0\leq s\leq m$ and $Q\in
C([0,m])$ by
\begin{equation}
    \pi(s,Q)\eqdef{{\gamma(0,Q(0))\over {\gamma(s,Q(s))}}}\,\exp\left\{-\int_0^s\frac{\mu(r,Q(r))}{\gamma(r,Q(r))}\,dr
    \right\}.
\end{equation}
\end{proposition}
\Pf{.} For a positive stationary solution $u_*$ let $Q_*$ be given
by \eqref{Qdef}. Since any stationary solution satisfies
\begin{align}\label{statsol2}
    u_*(s)=u_*(0)\,\pi(s,Q_*),
\end{align}
we obtain Eq.~\eqref{netrep} when imposing the boundary condition
\eqref{boundary}. \eofproof

Given any stationary solution $u_*$ in $C^1([0,m])$, we linearize
the governing equations by introducing the infinitesimal perturbation $v=v(s,t)$ and  making the ansatz
$u=v+u_*$. After inserting this expression in the governing equations and omitting all nonlinear terms,
we obtain the linearized problem
\begin{align}
 &v_t(s,t)+\gamma_*(s)\,v_s(s,t)+\rho_*(s)\,v(s,t) +\sigma_*(s)\,
 V(s,t)=0,\label{linear}\\
 &v(0,t)=\int_0^m \beta(s,Q_*(s))\,v(s,t)\,ds+\int_0^m\beta_Q(s,Q_*(s))\,u_*(s)\,V(s,t)\,ds,\label{linearboundary}
\end{align}
where we have set
\begin{align}
V(s,t)=&\,\alpha\,\int_0^s w(\eta)\,v(\eta,t)\,d\eta+\int_s^m
w(\eta)\,v(\eta,t)\,d\eta,\\
\gamma_*(s)=&\,\gamma(s,Q_*(s)),\\
\rho_*(s) =&\, \mu(s,Q_*(s))+\gamma_s(s,Q_*(s))+2\,(\alpha-1)\,w(s)\,\gamma_Q(s,Q_*(s))\,u_*(s),\\
\sigma_*(s)=& \,\mu_Q(s,Q_*(s))\,u_*(s)+\gamma_{sQ}(s,Q_*(s))\,u_*(s)+\gamma_Q(s,Q_*(s))\,u_*^\prime(s)\nonumber\\
&\quad +(\alpha-1)\,w(s)\,\gamma_{QQ}(s,Q_*(s))\,u_*(s)^2.
\end{align}
We denote the Lebesgue space $L^1(0,m)$ with its usual
norm $\|\cdot\|$ by ${\mathcal X}$ and introduce the bounded linear
functional $\Lambda$ on $\mathcal X$ by
\begin{align}
    \Lambda(v)= &\int_0^m \beta(s,Q_*(s))\,v(s)\,ds\label{Lambdadef}\\
    &+\int_0^m
    \beta_Q(s,Q_*(s))\,u_*(s)\,\left(\alpha\,\int_0^s w(\eta)\,v(\eta)\,d\eta+\int_s^m
w(\eta)\,v(\eta)\,d\eta\right)\,ds.\nonumber
\end{align}
Next we define the operators
\begin{align}
    {\mathcal A} v =& -\gamma(\cdot,Q_*)\,v_s,\quad
    \text{Dom}({\mathcal A}) =
    \left\{v\in W^{1,1}(0,m)\,|\,v(0)=\Lambda(v)\right\},\\
     {\mathcal B} v = &-\rho_*\,v\quad
    \text{on ${\mathcal X}$,}\label{defB}\\
     {\mathcal C} v= & -\sigma_*\,\left(\alpha\,\int_0^\cdot w(\eta)\,v(\eta)\,d\eta+\int_\cdot^m
    w(\eta)\,v(\eta)\,d\eta\right)\quad \text{on ${\mathcal X}$.}\label{defC}
\end{align}
Then the linearized system \eqref{linear}--\eqref{linearboundary} can be cast
in the form of an initial value problem for an
ordinary differential equation on ${\mathcal X}$
\begin{equation}\label{abstr}
    {d\over dt}\, v = \left({\mathcal A} + {\mathcal B}+ {\mathcal C}\right)\,v,
\end{equation}
together with the initial condition
\begin{equation}\label{abstr2}
    v(0)=v_0.
\end{equation}
In analogy to previously discussed size-structured population
models \cite{FH1,FH3,FH2}, we can invoke the
 Desch-Schappacher Perturbation Theorem \cite{NAG} to obtain the following result.

\begin{proposition}\label{gen}
The operator ${\mathcal A} + {\mathcal B}+ {\mathcal C}$ generates
a strongly continuous semigroup $\{{\mathcal T}(t)\}_{t\geq 0}$ of
bounded linear operators on ${\mathcal X}$.
\end{proposition}
The proof is a minor modification of parallel results given in
\cite{FH1,FH3,FH2} and has
therefore been omitted.

\sector{Spectral analysis and semigroup regularity}

\begin{proposition}\label{eigval}
    The spectrum of
    ${\mathcal A} + {\mathcal B}+ {\mathcal C}$ can contain only
    isolated eigenvalues of finite multiplicity.
\end{proposition}
\Pf{.} We  prove that the resolvent operator of ${\mathcal A} + {\mathcal B}+ {\mathcal
    C}$ is compact. Since ${\mathcal B}+ {\mathcal
    C}$ is a bounded perturbation of ${\mathcal A}$, it suffices to show that the resolvent operator of
    ${\mathcal A}$ is compact. To this end, given $f\in \mathcal X$, we find a unique solution $v\in \text{Dom}({\mathcal A})$
    of the equation
    \begin{equation}
        \lambda\,v-{\mathcal A}\,v= f
    \end{equation}
    in the form
    \begin{equation}
        v(s) =e^{-\lambda\,\Gamma(s)}\,\left(\Lambda(v)+\int_0^s e^{\lambda\,\Gamma(r)}\,{{f(r)}\over {
        \gamma(r,Q_*(r))}}\,dr\right)
    \end{equation}
    if $\lambda\in \mathbb R$ is sufficiently large. Here we define
    \begin{equation}
       \Gamma(s)\eqdef \int_0^s {1\over {\gamma_*(\eta)}}\,d\eta.
    \end{equation}
    Consequently, for $\lambda>0$ large enough,
     the resolvent operator $\left(\lambda\,{\mathcal I}-{\mathcal A}\right)^{-1}$ exists and is  bounded,
     mapping ${\mathcal X}=L^1(0,m)$ into
$W^{1,1}(0,m)$. Since
$W^{1,1}(0,m)$ is compactly embedded in $L^1(0,m)$, the claim
follows. \eofproof

\begin{theorem}\label{reg1}
    The semigroup $\{{\mathcal T}(t)\}_{t\geq 0}$, generated by the operator
    ${\mathcal A} + {\mathcal B}+ {\mathcal C}$, is eventually compact.
    Consequently, the Spectral Mapping Theorem holds true, i.e.
    \begin{equation}
        \sigma\left({\mathcal T}(t)\right) = \{0\}\cup \exp\{\sigma({\mathcal A} + {\mathcal B}+ {\mathcal C})\,t\},\quad t>0.
    \end{equation}
    Moreover, the semigroup is spectrally determined, i.e.\,the growth rate $\omega({\mathcal T})$ of the
    semigroup and the spectral bound $s({\mathcal A} + {\mathcal B}+ {\mathcal
    C})$ of its generator coincide.
\end{theorem}
\Pf{.}
Since the operator ${\mathcal C}$ is compact, it is enough to prove the claim for the operator
${\mathcal A}+{\mathcal B}$. The
differential equation
\begin{equation}\label{eqred}
    {d\over {dt}}\, v = ({\mathcal A}+{\mathcal B})\,v
\end{equation}
corresponds to the partial differential equation
\begin{equation}\label{reduc}
    v_t(s,t)+ \gamma_*(s)\,v_s(s,t)
    +\rho_*(s)\,v(s,t) = 0
\end{equation}
together with the boundary condition \eqref{linearboundary}.
For $t_0>0$ let us introduce
\begin{equation}
    \omega(s) = v(s,t(s)),
\end{equation}
where
\begin{equation}
    t(s)=t_0+ \Gamma(s).
\end{equation}
Then $\omega$ satisfies the equation
\begin{equation}
    \omega^\prime(s)+{{\rho_*(s)}\over {\gamma_*(s)}}\,\omega(s) = 0,
\end{equation}
hence
\begin{equation}
    \omega(s) = \Lambda(v(\cdot,t_0))\,\pi(s,Q_*)\,\exp\left((1-\alpha)\,\int_0^s
{{w(s)\,\gamma_Q(\eta,Q_*(\eta))\,u_*(\eta)} \over {\gamma_*(\eta)}}\,d\eta\right).
\end{equation}
Thus for $t-\Gamma(s) >0$ we have
\begin{equation}\label{vsmooth}
    v(s,t)= \Lambda(v(\cdot,t-\Gamma(s)))\,\pi(s,Q_*)\,\exp\left((1-\alpha)\,\int_0^s
{{w(s)\,\gamma_Q(\eta,Q_*(\eta))\,u_*(\eta)} \over {\gamma_*(\eta)}}\,d\eta\right).
\end{equation}
Therefore, noting the definition of $\Lambda$ in
\eqref{Lambdadef}, we conclude that $v$ is continuous in $s$ and
$t$ if $t>\Gamma(m) = \max_{0\leq s\leq m} \Gamma(s)$.
Consequently, Eq.~\eqref{vsmooth} in combination with Eq.~\eqref{reduc} implies that $v$ is continuously
differentiable if $t>2\,\Gamma(m)$. Hence the semigroup generated
by ${\mathcal A}+{\mathcal B}$ is differentiable for
$t>2\,\Gamma(m)$. Finally, since $W^{1,1}(0,m)$ is compactly
embedded in $L^1(0,m)$, the semigroup is compact for
$t>2\,\Gamma(m)$. The validity of the Spectral Mapping Theorem and
the claim about the spectral determinacy of the semigroup  follow, see
\cite{NAG}. \eofproof

We conclude this section by formulating conditions for the
positivity of the semigroup $\{{\mathcal T}(t)\}_{t\geq 0}$.
\begin{theorem}\label{pos}
Suppose that
\begin{align}
    &\sigma_*\leq 0\quad \text{and}\label{pos1} \\
    & \beta(\cdot, Q_*)+w\,\left(\int_0^\cdot \beta_Q(\eta,Q_*(\eta))\,u_*(\eta)\,d\eta +
    \alpha\,\int_\cdot^m\beta_Q(\eta,Q_*(\eta))\,u_*(\eta)\,d\eta
    \right)\geq 0.\label{pos2}
\end{align}
Then the semigroup $\{{\mathcal T}(t)\}_{t\geq 0}$, generated by
the operator ${\mathcal A} + {\mathcal B}+ {\mathcal C}$, is
positive.
\end{theorem}
\begin{remark}
  Conditions~\eqref{pos1}, \eqref{pos2} are immediate generalizations of the corresponding positivity conditions
  given by Pr\"u\ss\ in \cite{PR1} for an age-structured scramble competition model. In general, if $\beta_Q\equiv 0$, condition \eqref{pos2}
  is trivially satisfied. Also, if the growth rate is independent of the environment (i.e. $\gamma=\gamma(s)$), condition \eqref{pos1} reduces to
 \begin{equation}
   \mu_Q(s,Q_*(s))\le 0,\quad s\in [0,\infty).
\end{equation}
Hence in this case mortality is required to be a non-increasing
function of the environment as well.
\end{remark}

\Pf{\ of Theorem~\ref{pos}}
Since $\mathcal C$ is a positive operator by
condition \eqref{pos1}, it suffices to prove the claim for the semigroup generated by
${\mathcal A}+{\mathcal B}$. Hence we assume that $v$ satisfies Eq.~\eqref{reduc} such that
the boundary condition \eqref{linearboundary} and the initial condition
$v=v_0\in \text{Dom}({\mathcal A})$ for $t=0$ hold true.
Let the function $e_*$ be given by
\begin{equation}
    e_*(s) = \exp\left((1-\alpha)\,\int_0^s
{{w(s)\,\gamma_Q(\eta,Q_*(\eta))\,u_*(\eta)} \over {\gamma_*(\eta)}}\,d\eta\right).
\end{equation}
Then the function $\phi$, defined by
\begin{equation}
        \phi(s,t)={{v(s,t)}\over {\pi(s,Q_*)\,e_*(s)}},
\end{equation}
solves the problem
\begin{align}
    & \phi_t(s,t)+ \gamma(s,Q_*(s))\,\phi_s(s,t) = 0,\label{ueq1}\\
    & \phi(0,t)=\Lambda\left({{\phi(\cdot, t)}}\,\pi(\cdot,Q_*)\,e_*\right),\label{ueq2}\\
    & \phi(s,0) = {{v_0(s)}\over {\pi(s,Q_*)\,e_*(s)}}\eqdef \phi_0(s).\label{ueq3}
\end{align}
This boundary-initial value problem corresponds to the abstract
initial value problem
\begin{equation}
        {d\over dt}\, \phi = {\mathcal A}_M\,\phi,\quad \phi(0)=\phi_0
\end{equation}
with the modified semigroup generator ${\mathcal A}_M$, defined by
\begin{align}
    {\mathcal A}_M \phi =& -\gamma(\cdot,Q_*)\,\phi_s\quad \text{on the domain}\nonumber\\
    &\text{Dom} ({\mathcal A}_M) =\left\{
    \phi\in W^{1,1}(0,m)\,|\, \phi(0)=\Lambda\left({\phi}\,\pi(\cdot,Q_*)\,e_*\right)\right\}.
\end{align}
For $\lambda\geq 0$ and $g\in L^1(0,m)$, the
resolvent equation
\begin{equation}
    \lambda \phi -{\mathcal A}_M \phi = g
\end{equation}
has the implicit solution
\begin{equation}\label{usol}
    \phi(s)=e^{-\lambda\,\Gamma(s)}\,\Lambda\left({\phi}\,\pi(\cdot,Q_*)\,e_*\right)+ \int_0^s e^{\lambda\,(\Gamma(r)-\Gamma(s))}\,{{g(r)}\over
    {\gamma(r,Q_*(r))}}\,dr.
\end{equation}
Applying $\Lambda$, we deduce the equation
\begin{equation}
    \Lambda\left({\phi}\,\pi(\cdot,Q_*)\,e_*\right) =
    {\displaystyle{\Lambda\left(\int_0^\cdot e^{\lambda\,(\Gamma(r)-\Gamma(\cdot))}\,{{g(r)}\over
    {\gamma(r,Q_*(r))}}\,dr\,{\pi(\cdot,Q_*)\,e_*}\right)}\over
    {1-\Lambda\left({{e^{-\lambda\,\Gamma}}\,\pi(\cdot,Q_*)\,e_*}\right)}}
\end{equation}
if $\lambda$ is large enough. Condition \eqref{pos2} guarantees that $\Lambda$ is
a positive linear functional. Hence the solution $\phi$, given by
Eq.~\eqref{usol}, is nonnegative if $g$ is nonnegative and $\lambda$ is sufficiently
large. It follows that the resolvent operator
of ${\mathcal A}_M$ (and consequently of ${\mathcal A}+{\mathcal
B}$) is positive if $\lambda$ is large enough. This observation proves the claim.
\eofproof

The positivity of the semigroup has far-reaching consequences. In
particular, we obtain the following result from the theory of
positive semigroups \cite{NAG}.

\begin{corollary}\label{poscor}
    Suppose that the semigroup $\{{\mathcal T}(t)\}_{t\geq 0}$, generated by the operator ${\mathcal A} + {\mathcal B}+ {\mathcal C}$, is
    positive.
    Then the spectral bound $s({\mathcal A}+{\mathcal B}+{\mathcal C})\in [-\infty, \infty)$ satisfies
    \begin{equation}
        s({\mathcal A}+{\mathcal B}+{\mathcal C}) = \max\,\left\{\lambda\in {\mathbb R}\,|\,\text{$\lambda$ is
        eigenvalue of ${\mathcal A}+{\mathcal B}+{\mathcal C}$}\right\}.
    \end{equation}
    Moreover,  the spectrum of ${\mathcal A}+{\mathcal B}+{\mathcal C}$ is nonempty if and only if the spectral bound is finite.
\end{corollary}

\sector{The characteristic equation}

In light of Theorem~\ref{reg1} the growth of the governing
semigroup is determined by the eigenvalues of its generator.
Hence
it is essential to determine the eigenvalues of ${\mathcal A} +
{\mathcal B}+ {\mathcal C}$. The eigenvalue equation
\begin{equation}
    \lambda\,v-\left({\mathcal A} + {\mathcal B}+ {\mathcal C}\right) v =0
\end{equation}
for $\lambda\in\mathbb C$ and nontrivial $v$ is equivalent to the system
\begin{align}
 &v^\prime(s)\,\gamma_*(s)+v(s)\,(\lambda+\rho_*(s))+V(s)\,\sigma_*(s)=0,\label{linode1}\\
 &v(0)=\int_0^m\beta(s,Q_*(s))\,v(s)\,ds+\int_0^m\beta_Q(s,Q_*(s))\,u_*(s)\,V(s)\,ds,\label{linode3}
\end{align}
where
\begin{align}
V(s)=&\ \alpha\int_0^s w(\eta)\,v(\eta)\,d\eta+\int_s^m w(\eta)\,v(\eta)\,d\eta \nonumber\\
=&\ (\alpha-1)\,\int_0^s w(\eta)\,v(\eta)\,d\eta+\int_0^m w(\eta)\,v(\eta)\,d\eta.\label{Veq}
\end{align}
For the remainder of this section let us assume that $\alpha\in [0,1)$.
From \eqref{Veq} we obtain
\begin{equation}\label{Veq2}
V^\prime(s)=(\alpha-1)\,w(s)\,v(s)\quad \text{and}\quad  V^{\prime\prime}(s)=(\alpha-1)\,(w^\prime(s)\,v(s)+w(s)\,v^\prime(s)).
\end{equation}
Using the relations \eqref{Veq2} we can rewrite system \eqref{linode1}--\eqref{linode3} in terms of $V$ and its derivatives as follows
\begin{equation}
  V^{\prime\prime}(s)+V^\prime(s)\,\left(\frac{\rho_*(s)+\lambda}
  {\gamma_*(s)}-\frac{w'(s)}{w(s)}\right)
  +V(s)\,(\alpha-1)\,w(s)\,{{\sigma_*(s)}\over {\gamma_*(s)}}=0.\label{transformed}
\end{equation}
Eq.~\eqref{transformed} is accompanied by boundary conditions of the form
\begin{align}
&\alpha V(0)=V(m),\label{transbound1}\\
&V^\prime(0)=w(0)\,\int_0^m\frac{\beta(s,Q_*(s))}{w(s)}\,V'(s)\,ds\nonumber\\
&\quad +w(0)\,(\alpha-1)\,\int_0^m\beta_Q(s,Q_*(s))\,u_*(s)\,V(s)\,ds.\label{transbound2}
\end{align}
For $\lambda\in \mathbb C$, any solution $V_\lambda(s)$ of
the second order homogeneous ordinary differential
equation
\eqref{transformed} can be written as
\begin{equation}
V_\lambda(s)=c_1\,V_1(s,\lambda)+c_2\,V_2(s,\lambda),
\end{equation}
where $V_1(s,\lambda)$ and $V_2(s,\lambda)$ are any fixed,
linearly independent solutions of Eq.~\eqref{transformed} and
$c_1$, $c_2$ are arbitrary constants. When imposing the boundary
conditions \eqref{transbound1}--\eqref{transbound2}, we obtain the
conditions
\begin{align}
 & c_1\,V^\prime_1(0,\lambda)+c_2\,V^\prime_2(0,\lambda)=c_1\,
 \int_0^m\frac{w(0)}{w(s)}\,\beta(s,Q_*(s))\,V^\prime_1(s,\lambda)\,ds\\
 &+c_2\,\int_0^m\frac{w(0)}{w(s)}\,\beta(s,Q_*(s))\,V^\prime_2(s,\lambda)\,ds\\
 &+c_1\,\int_0^m(\alpha-1)\,w(0)\,\beta_Q(s,Q_*(s))\,u_*(s)\,V_1(s,\lambda)\,ds\\
 &+c_2\,
 \int_0^m(\alpha-1)\,w(0)\,\beta_Q(s,Q_*(s))\,u_*(s)\,V_2(s,\lambda)\,ds
\end{align}
or in short
\begin{equation}
    c_1\,H_1(\lambda)+c_2\,H_2(\lambda) = 0,\label{h1}
\end{equation}
and
\begin{equation}
 c_1\,\alpha\, V_1(0,\lambda)+c_2\,\alpha\, V_2(0,\lambda)=c_1\,V_1(m,\lambda)+c_2\,V_2(m,\lambda),
\end{equation}
in short
\begin{equation}
    c_1\,J_1(\lambda)+c_2\,J_2(\lambda) = 0.\label{h2}
\end{equation}
Here the functions $H_1$, $H_2$, $J_1$, and $J_2$ represent the
terms multiplying $c_1$, $c_2$, respectively. The homogeneous
system \eqref{h1}, \eqref{h2}  admits a nontrivial solution for
$c_1$, $c_2$ if and only if $\lambda$ satisfies the equation
\begin{equation}
    H_1(\lambda)\,J_2(\lambda)-H_2(\lambda)\,J_1(\lambda)=0.\label{genchareq}
\end{equation}
This equation is the characteristic equation of the linearized
system \eqref{linear}-\eqref{linearboundary}. Its zeros are the
eigenvalues of the operator ${\mathcal A} + {\mathcal B}+
{\mathcal C}$, which completely describe the spectrum of
${\mathcal A} + {\mathcal B}+ {\mathcal C}$.

The explicit information contained in the characteristic equation
is, however, rather limited since linearly independent solutions
of the second order differential equation \eqref{transformed} are
in general not directly available, unless one resorts to numerical
techniques. As we will see in the forthcoming section this problem
can, however, be overcome in special cases of the model
ingredients.

\sector{A special case}

In this section we treat the special case when the mortality and growth rate are
independent of the environment $Q$, i.e.\,$\gamma_Q\equiv 0\equiv \mu_Q$. Hence $\sigma_*\equiv 0$ and
$e_*\equiv 1$.
 In this situation we are able to
determine the characteristic equation \eqref{genchareq} explicitly
and to formulate simple conditions for the linear
stability/instability of positive stationary solutions.  In contrast to the preceding section we allow $\alpha\in
[0,1]$.

\begin{theorem}\label{specstab}
Suppose $\sigma_*\equiv 0$. Then a positive stationary solution $u_*$
is linearly asymptotically stable if
\begin{equation}
\beta_Q(\cdot,Q_*)\leq 0,\quad \beta_Q(\cdot, Q_*)\not\equiv 0
\label{stabcond1}
\end{equation}
and the positivity condition \eqref{pos2} holds true. If, however,
\begin{equation}\label{instab}
\beta_Q(.,Q_*)\geq 0,\quad \beta_Q(\cdot, Q_*)\not\equiv 0,
\end{equation}
then $u_*$ is linearly unstable.
\end{theorem}
Note that the instability part of the theorem does not require the positivity condition.

\Pf{.} We assume first that $0\leq \alpha<1$. Then the general solution of
\eqref{transformed} is found as
\begin{equation}
V(s)=V(0)+V^\prime(0)\int_0^s\frac{w(r)}{w(0)}\,\Pi(\lambda,r)\,dr,\label{specsolution}
\end{equation}
where we have set
\begin{equation}
\Pi(\lambda,r)\eqdef{{\gamma_*(0)\over {\gamma_*(s)}}}\,\exp\left\{-\int_0^s
\frac{\lambda+\mu_*(r)}{\gamma_*(r)}\,dr \right\}.
\end{equation}
Imposing the boundary condition \eqref{transbound1} on the
solution \eqref{specsolution}, we obtain
\begin{equation}
0=V(0)\,(1-\alpha)+V^\prime(0)\,\int_0^m\frac{w(s)}{w(0)}\,\Pi(\lambda,s)\,ds,\label{homog1}
\end{equation}
while the boundary condition \eqref{transbound2} gives
\begin{align}
 0=&\ V(0)\,\left(w(0)\,(1-\alpha)\int_0^m\beta_Q(s,Q_*(s))\,u_*(s)\,ds\right)\nonumber\\
&\ +V^\prime(0)\,\left(1-\int_0^m\beta(s,Q_*(s))\,\Pi(\lambda,s)\,ds\right.\label{homog2}\\
&\ \left.+(1-\alpha)\,\int_0^m\beta_Q(s,Q_*(s))\,u_*(s)\,\int_0^s
w(r)\,\Pi(\lambda,r)\,dr\,ds\right).\nonumber
\end{align}
The linear system \eqref{homog1}--\eqref{homog2} has a nontrivial
solution $(V(0),V'(0))$ if and only if $\lambda$ satisfies
\begin{align}
1= & \int_0^m\beta(s,Q_*(s))\,\Pi(\lambda,s)\,ds\nonumber\\
& +(\alpha-1)\,\int_0^m\beta_Q(s,Q_*(s))\,u_*(s)\,\int_0^s w(r)\,\Pi(\lambda,r)\,dr\,ds\label{specchareq}\\
& +\int_0^m\beta_Q(s,Q_*(s))\,u_*(s)\,ds\,\int_0^m
w(s)\,\Pi(\lambda,s)\,ds\eqdef K(\lambda).\nonumber
\end{align}
This equation corresponds to the characteristic equation
\eqref{genchareq}.
If, however, $\alpha=1$, $V$, defined by \eqref{Veq}, is constant. When we solve the problem \eqref{linode1}--\eqref{linode3}
directly, we obtain again the condition
\begin{equation}
    K(\lambda)=1,
\end{equation}
where $K$ is given by \eqref{specchareq} with $\alpha=1$. Hence \eqref{specchareq}
is the characteristic equation for all $0\leq \alpha\leq 1$.
For the stability part, our assumptions guarantee that the
positivity conditions \eqref{pos1}, \eqref{pos2} hold true.
Therefore, to prove asymptotic stability, it suffices to show that
the characteristic equation \eqref{specchareq} has  no nonnegative
(real) solutions. To this end, we observe that
\begin{align}
K(0)= & R(Q_*)
+\int_0^m\beta_Q(s,Q_*(s))\,u_*(s)\nonumber\\
&\times\left(\alpha\,\int_0^s w(r)\,\pi(r,Q_*)\,dr+\int_s^m
w(r)\,\pi(r,Q_*)\,dr\right)\,ds<1
\end{align}
by condition \eqref{stabcond1}. Moreover, the positivity condition
\eqref{pos2} yields that
\begin{align}
K^\prime(\lambda)= & -\int_0^m\Pi(\lambda,s)\,\int_0^s\frac{1}{\gamma(r,Q_*(r))}\,dr\,\bigg(\beta(s,Q_*(s))\nonumber\\
&+w(s)\,\int_0^s\beta_Q(r,Q_*(r))\,u_*(r)\,dr \label{Kder} \\
&\left.+\alpha\,w(s)\int_s^m\beta_Q(r,Q_*(r))\,u_*(r)\,dr\right)\,ds\leq
0. \nonumber
\end{align}
Consequently, $K$ is monotone decreasing for $\lambda\ge 0$. Hence
the stability part is proven. The instability part of
the theorem follows from the Intermediate Value Theorem since
$K(0)>1$ by \eqref{instab} and  $\displaystyle\lim_{\lambda\rightarrow \infty} K(\lambda)=0$.
\eofproof

\begin{example}
Let us consider an example where Theorem~\ref{specstab} yields
asymptotic stability. We choose
\begin{equation}
    m=1,\quad \alpha={1\over 2},\quad w\equiv 1
\end{equation}
and let
\begin{align}
    &\gamma(s) = 1-{1\over 2}\,s,\qquad \mu\equiv 1,\\
    &\beta(s,Q) = {480\over 997}\,(1+s)\,(3-2\,Q)\quad \text{if
    $Q\leq {3\over 4}$, $0\leq s\leq 1$,}
\end{align}
where we assume that $\beta$ extends to a continuously
differentiable, non-negative function on $[0,1]\times [0,\infty)$.
Then the corresponding problem has the stationary solution
$u_*(s)=1-{1\over 2}\,s$ with
\begin{equation}
    Q_*(s) = {{s^2}\over 8}-{s\over 2}+{3\over 4}\leq {3\over 4}\quad \text{for $0\leq s\leq
    1$.}
\end{equation}
It is readily seen that
\begin{equation}
    \beta_Q(s,Q_*(s)) = -{960\over 997}\,(1+s)<0,\quad 0\leq s\leq
    1
\end{equation}
and that the positivity condition \eqref{pos2} reduces to
\begin{equation}
    -{{s^3}\over 24}+{{s^2}\over 4}+{{3\,s}\over
    4}+{5\over 24}\geq 0,\quad 0\leq s\leq 1.
\end{equation}
Since this inequality holds true, the stationary solution $u_*$ is
linearly asymptotically stable.
\end{example}

\sector{Direct approach: dissipativity}

Since in the general case of environment dependent vital rates  the characteristic
equation is not explicitly available, we shall pursue a different
path to obtain asymptotic stability. This approach is based on
dissipativity calculations in the underlying state space
${\mathcal X}=L^1(0,m)$ and proceeds parallel to
similar developments for the linear semigroup of fiber spinning in
\cite{Ha2,Ha4}. An added advantage of this
technique is that we can discuss linear stability of the trivial stationary solution
and that we can forego imposing positivity conditions on
the semigroup. In addition we can include the case $\alpha=1$
without technical difficulties. To our knowledge dissipativity
estimates have so far not been used in the case of hierarchical
size-structured population models.

\begin{theorem}\label{estimate}
A stationary solution $u_*$ is linearly asymptotically stable if
\begin{align}
    &\mu(s,Q_*(s))> w(s)\,\left((1-\alpha)\,\gamma_Q(s, Q_*(s))\,u_*(s)+||\sigma_*||\big)\right)\nonumber\\
    &\quad +
    \gamma(0,Q_*(0))\,\bigg|\beta(s,Q_*(s))
    +\alpha\,w(s)\int_s^m\beta_Q(r,Q_*(r))\,u_*(r)\,dr \label{disscond}\\
    &\quad +w(s)\,\int_0^s\beta_Q(r,Q_*(r))\,u_*(r)\,dr\bigg|\nonumber
\end{align}
for $0\leq s\leq m$.
\end{theorem}
As before the norm on $L^1(0,m)$ is denoted by $\|\cdot\|$.

\Pf{.}
We will show that, under the given condition, there exists $\kappa>0$ such
that the operator ${\mathcal A} +{\mathcal B}+ {\mathcal C}+\kappa\,{\mathcal I}$ is dissipative.
Consequently,  the semigroup $\{{\mathcal T}(t)\}_{t\geq 0}$ generated by the operator
    ${\mathcal A} + {\mathcal B}+ {\mathcal C}$ obeys
\begin{equation}
    \|{\mathcal T}(t)\|\leq e^{-\kappa\,t},\quad t\geq 0,
\end{equation}
which proves the claim.

To obtain dissipativity, assume that, for given $h\in \mathcal X$, $v\in \text{Dom}({\mathcal A})$ is such that, for some $\lambda>0$,
\begin{equation}
    v-\lambda\,({\mathcal A} +{\mathcal B}+ {\mathcal C}+\kappa\,{\mathcal I})\,v=h.\label{res}
\end{equation}
Then we have
\begin{align}
    \|v\| =&\ \int_0^m v(s)\,\text{sgn}\,v(s)\,ds\nonumber\\
    =&\ \int_0^m h(s)\,\text{sgn}\,v(s)\,ds-\lambda\,\int_0^m \gamma(s,Q_*(s))\,v_s(s)\,\text{sgn}\,v(s)\,ds\nonumber\\
    &\ -\lambda\,\int_0^m \rho_*(s)\,v(s)\,\text{sgn}\,v(s)\,ds\nonumber\\
    &\ -\lambda\,\int_0^m \sigma_*(s)\,\left(\alpha\,\int_0^s w(\eta)\,v(\eta)\,d\eta\right.\\
    &\ \left.+\int_s^m
    w(\eta)\,v(\eta)\,d\eta\right)\,\text{sgn}\,v(s)\,ds\nonumber+\lambda\,\kappa\,\int_0^m v(s)\,\text{sgn}\,v(s)\,ds.
\end{align}
Here we have used the definition $\text{sgn}\,0 = 0$.
The set of points in the interval $(0,m)$ where $v$ is nonzero is the countable union of disjoint open intervals $(a_i,b_i)$ on each of which either $v>0$ or $v<0$
holds true such that $v(a_i)=0$ for all $i$ unless $a_i=0$, and such that $v(b_i)=0$ unless $b_i=m$.
If $(a_i,b_i)$ is any such interval on which $v>0$ we have after integration by parts
\begin{align}
    &\int_{a_i}^{b_i} v(s)\,ds\leq
    \int_{a_i}^{b_i}|h(s)|\,ds-\lambda\,\gamma(b_i,Q_*(b_i))\,v(b_i)+\lambda\,\gamma(a_i,Q_*(a_i))\,v(a_i)\nonumber \\
   &\quad+\lambda\,\int_{a_i}^{b_i}\left(\kappa-\mu(s,Q_*(s))-(\alpha-1)\,w(s)\,\gamma_Q(s, Q_*(s))\,u_*(s)\right)\,v(s)\,ds\\
    &\quad +
    \lambda\,\int_{a_i}^{b_i}|\sigma_*(s)|\,ds\,\int_0^m w(s)\,|v(s)|\,ds.\nonumber
\end{align}
Similarly, on any interval $(a_i,b_i)$ where $v<0$ we have
\begin{align}
& \int_{a_i}^{b_i} |v(s)|\, ds\leq \int_{a_i}^{b_i}|h(s)|\,ds+
\lambda\,\gamma(b_i,Q_*(b_i))\,v(b_i)-\lambda\,\gamma(a_i,Q_*(a_i))\,v(a_i)\nonumber\\
&\quad+\lambda\,\int_{a_i}^{b_i}\left(\kappa-\mu(s,Q_*(s))-(\alpha-1)\,w(s)\,\gamma_Q(s, Q_*(s))\,u_*(s)\right)\,|v(s)|\,ds
\\
&\quad +\lambda\,\int_{a_i}^{b_i}|\sigma_*(s)|\,ds\,\int_0^m w(s)\,|v(s)|\,ds.\nonumber
\end{align}
Finally, noting that $v(a_i)=0=v(b_j)$ unless $a_i=0$, $b_j=m$, we combine these two estimates to obtain
\begin{align}
||v||\leq& \ ||h||+ \lambda\,\gamma(0,Q_*(0))\,|v(0)|\\
&\hspace*{-1cm} +\lambda\,\int_0^m\big(\kappa-\mu(s,Q_*(s))+w(s)\,\left((1-\alpha)\,\gamma_Q(s, Q_*(s))\,u_*(s)+||\sigma_*||\big)\right)\,|v(s)|\,ds.\nonumber
\end{align}
Since
\begin{align}
    |v(0)| =  |\Lambda(v)|\leq &\ \int_0^m \left|\beta(s,Q_*(s))
    +\alpha\,w(s)\int_0^m\beta_Q(r,Q_*(r))\,u_*(r)\,dr\right.\nonumber\\
    &\ \left.+(1-\alpha)\,w(s)\,\int_0^s\beta_Q(r,Q_*(r))\,u_*(r)\,dr\right|\,|v(s)|\,ds
\end{align}
and since condition \eqref{disscond} is satisfied, we can choose
$\kappa>0$ such that, for $0\leq s\leq m$,
\begin{align}
    &\kappa-\mu(s,Q_*(s))+w(s)\,\left((1-\alpha)\,\gamma_Q(s, Q_*(s))\,u_*(s)+||\sigma_*||\big)\right)\nonumber\\
    &\quad +
    \gamma(0,Q_*(0))\,\bigg|\beta(s,Q_*(s))+\alpha\,w(s)\int_0^m\beta_Q(r,Q_*(r))\,u_*(r)\,dr\\
    &\quad +(1-\alpha)\,w(s)\,\int_0^s\beta_Q(r,Q_*(r))\,u_*(r)\,dr\bigg|\leq 0.\nonumber
\end{align}
For such $\kappa$, we have the desired inequality
\begin{equation}
    \|v\|\leq \|h\|\quad\text{for $\lambda>0$,}
\end{equation}
thus establishing dissipativity. Hence we conclude that ${\mathcal
A} +{\mathcal B}+ {\mathcal C}+\kappa\,{\mathcal I}$ generates a
contraction semigroup.
\eofproof\\

\begin{remark}
    Suppose $u_*$ is a stationary solution such that condition
    \eqref{disscond} is satisfied with $\beta_Q\geq 0$.
   Then we have
   \begin{align}
        R(Q_*) = &\ \int_0^m
        {{\beta(s,Q_*(s))\,\gamma(0,Q_*(0))}\over
        {\gamma(s,Q_*(s))}}\,\exp\left(-\int_0^s
        {{\mu(r,Q_*(r))}\over
        {\gamma(r,Q_*(r))}}\,dr\right)\,ds\nonumber\\
        \leq&\ \int_0^m {{\mu(s,Q_*(s))}\over
        {\gamma(s,Q_*(s))}}\,\exp\left(-\int_0^s
        {{\mu(r,Q_*(r))}\over
        {\gamma(r,Q_*(r))}}\,dr\right)\,ds\\
        =&\ 1-\exp\left(-\int_0^m
        {{\mu(r,Q_*(r))}\over
        {\gamma(r,Q_*(r))}}\,dr\right)<1.
\end{align}
    Hence in light of \eqref{netrep} we have to conclude that
    $u_*\equiv 0$.
\end{remark}

\begin{remark}
For the stability of the trivial equilibrium $u_*\equiv 0$ the criterion \eqref{disscond} reduces to
\begin{equation}
\mu(s,0)>\gamma(0,0)\,\beta(s,0),\quad s\in [0,m].\label{trivial}
\end{equation}
 Note that
 \eqref{trivial} clearly implies $R(0)<1$, which is the well-known stability criterion of the trivial steady state
in scramble competition, see \cite{I}.
\end{remark}

\begin{remark}
In scramble competition ($\alpha=1$) the stability criterion
\eqref{disscond} for a stationary solution $u_*$ with
total (weighted) population
\begin{equation}
    P_*   =\int_0^m w(\eta)\,u_*(\eta)\,d\eta
\end{equation}
reads
\begin{equation}
\mu(s,P_*) >  \ w(s)\,||\sigma_*||+
    \bigg|\tilde{\beta}(s,P_*)
    +w(s)\,\int_0^m\tilde{\beta}_P(r,P_*)\,u_*(r)\,dr\bigg|,\ 0\leq s\leq m,\label{disscondred}
\end{equation}
where
\begin{equation}
    \tilde{\beta}(s,P) = \gamma(0,P)\,\beta(s,P).
\end{equation}
\end{remark}

\begin{example}
    We will give a nontrivial example of a stationary solution
    for which the stability criterion \eqref{disscond} holds true.
    We choose
\begin{equation}
    m=1,\quad \alpha={1\over 2},\quad w\equiv 1,\quad \gamma(s) = 1-{1\over 2}\,s,\quad \mu\equiv 1
\end{equation}
and let $\beta\in C^1([0,1]\times [0,\infty))$ be positive such that
\begin{equation}
    \beta(s,Q) = {160\over 159}\,(1+s)\,(2-2\,Q)\quad \text{if
    $Q\leq {3\over 4}$, $0\leq s\leq 1$.}
\end{equation}
Then we have again the stationary solution
$u_*(s)=1-{1\over 2}\,s$ with
\begin{equation}
    Q_*(s) = {{s^2}\over 8}-{s\over 2}+{3\over 4}\leq {3\over 4}\quad \text{for $0\leq s\leq
    1$.}
\end{equation}
Now, however, the positivity condition \eqref{pos2} is violated. Nonetheless we obtain
\begin{align}
    & \bigg|\beta(s,Q_*(s))+{1\over 2}\,\int_s^1\beta_Q(r,Q_*(r))\,u_*(r)\,dr +\int_0^s\beta_Q(r,Q_*(r))\,u_*(r)\,dr\bigg|\nonumber\\
    &\quad = {160\over 159}\,\bigg|{{s^3}\over {12}}-{{s^2}\over 2}-{s\over 2}+{7\over 12}\bigg|<1,
    \qquad 0\leq s\leq 1.
\end{align}
Hence the stationary solution is linearly asymptotically stable by Theorem~\ref{estimate}.
A straightforward perturbation argument can be used to extend this example to a more complicated
situation with environment dependent mortality and growth rate.
\end{example}

\sector{Conclusion}

In this work we have analyzed the linear asymptotic stability of
equilibrium solutions of a nonlinear hierarchical size-structured
population model. We have extended our previous mathematical
approach in \cite{FH1,FH3,FH2} for the case of scramble
competition models to the hierarchical case. As conjectured in
\cite{DGM} for general physiologically structured models, we note
that the linear asymptotic stability of stationary solutions is
determined by zeros of a characteristic function. When the linear
dynamical behavior is governed by a positive semigroup, the
characteristic function has a dominant real root, unless the
spectrum of the semigroup generator is empty. As we have seen in
Section 4, however, this function is not explicitly available
(except in special cases). Nevertheless, we managed to
characterize the spectrum of the linearized operator implicitly,
by deducing an eigenvalue problem for a second order differential
operator. This characterization allows in principle to further
investigate stability questions by numerical techniques in case of
concrete model ingredients.

To overcome the severe limitations caused by the spectral
characterizations of asymptotic stability, we have given a direct
dissipativity condition in terms of the model ingredients in the
relevant state space, guaranteeing the exponential decay of the
governing linear semigroup. This elementary, though important
criterion allows us to expand linear stability studies beyond the
setting of positive semigroups and dominant eigenvalues.

\bigskip

{\bf Acknowledgment}

JZF was supported by EPSRC grant EP/F025599/1. TH acknowledges
support through NSF-Grant DMS 0709197.

\end{document}